\documentclass[%a4paper,
11pt,oneside]{article}
\usepackage[english]{babel}

\usepackage[curve,matrix,arrow,cmtip]{xy}
\NoComputerModernTips

\xyoption{rotate}

\newdir^{ (}{{}*!/-3pt/\dir^{(}}    
\newdir_{ (}{{}*!/-3pt/\dir_{(}}

\usepackage[ansinew]{inputenc}
\usepackage[dvipsnames,usenames]{color} 

\usepackage[backend=bibtex, maxbibnames=50, style=alphabetic]{biblatex}
\long\def\forget#1{}%{{\rm #1}}

\usepackage[scale=0.7,hmarginratio=1:1,vmarginratio=2:3]{geometry}

\newcounter{themargin}
\usepackage{lmodern} 
\usepackage{diagbox}
\usepackage{latexsym}
\usepackage{amsmath}
\usepackage{amstext}
\usepackage{amssymb}
\usepackage{amsbsy}

\usepackage{amsthm}

\usepackage[active]{srcltx}
\usepackage[titles]{tocloft}

\usepackage{GBplainAutF}

\usepackage{enumitem}
\newbox\dottobox
\setbox\dottobox=\hbox{$\longrightarrow$}
\setbox\dottobox=\hbox to\wd\dottobox{\hss$
\UseComputerModernTips\xymatrix@C=.5cm{\ar@{.>}[r]&\\}
                                      $\hss}

\newbox\leftdottobox
\setbox\leftdottobox=\hbox{$\longleftarrow$}
\setbox\leftdottobox=\hbox to\wd\leftdottobox{\hss$
\UseComputerModernTips\xymatrix@C=.5cm{\ar@{<.}[r]&\\}
                                      $\hss}

\newbox\dotintobox
\setbox\dotintobox=\hbox{$\longrightarrow$}
\setbox\dotintobox=\hbox to\wd\dotintobox{\hss$
\UseComputerModernTips\xymatrix@C=.5cm{\ar@{^{ (}.>}[r]&\\}
                                      $\hss}

\usepackage{lmodern,graphicx}

\usepackage{geboeautoref-section}

\theoremstyle{plain} 
\maketheorem{ThmDef}{Theorem-Definition}{Thm}
\maketheorem{Prop}{Proposition}{Thm}
\maketheorem{PropDef}{Proposition--Definition}{Thm}
\maketheorem{LemDef}{Lemma--Definition}{Thm}
\maketheorem{KeyLem}{Key Lemma}{Thm}
\maketheorem{Clm}{Claim}{Thm}
\maketheorem{Lem}{Lemma}{Thm}
\maketheorem{MainLem}{Main lemma}{Thm}
\maketheorem{SubLem}{Sublemma}{Thm}
\maketheorem{Cor}{Corollary}{Thm}
\maketheorem{Fact}{Fact}{Thm}
\maketheorem{Conj}{Conjecture}{Thm}
\maketheorem{Def}{Definition}{Thm}
\maketheorem{Defs}{Definitions}{Thm}

\newtheorem*{ThmInt}{Theorem}

\theoremstyle{definition} 
\maketheorem{Conv}{Convention}{Thm}
\maketheorem{Alg}{Algorithm}{Thm}
\maketheorem{Ex}{Example}{Thm}
\maketheorem{Exs}{Examples}{Thm}

\theoremstyle{remark} 
\maketheorem{Rem}{Remark}{Thm}
\maketheorem{Rems}{Remarks}{Thm}
\maketheorem{Parag}{}{Thm}
\maketheorem{Empt}{}{Thm}
\maketheorem{Exer}{Exercise}{Thm}
\maketheorem{Prob}{Problem}{Thm}
\maketheorem{Cond}{Condition}{Thm}
\maketheorem{Conds}{Conditions}{Thm}
\maketheorem{Ass}{Assumption}{Thm}
\maketheorem{Que}{Question}{Thm}
\usepackage{comment}
\usepackage{csquotes}
\newcommand{\har}{\mathrm{har}}

\newcommand{\Tempnewpage}%{}%
{\newpage}

\newcommand{\ro}{\varphi}
\newcommand{\ri}{\widehat{\varphi}}
\newcommand{\ol}[1]{{\overline{#1}}}

\title{Dimension formulas for certain spaces of Drinfeld cusp forms}
\author{Gebhard\ B\"ockle, Peter Mathias Gr\"af and Iason Papadopoulos}
\date{\today}
\begin{document}
\maketitle

\setcounter{secnumdepth}{3}
\setcounter{tocdepth}{3}

\setcounter{section}{0}

\begin{abstract}
In this short note, we derive dimension formulas for spaces of Drinfeld cusp forms corresponding to harmonic cocycles invariant under the group $\SL_2(\BF_q[t])$ and with values in absolutely irreducible $\SL_2(\BF_q(t))$-representations via the theory of Brauer characters. This generalizes results in \cite{BGP} obtained by different methods. In addition, we prove a simple asymptotic formula for these dimensions.
\end{abstract}

\section{Introduction}

Let $F$ be a global function field with finite field of constants $\BF_q$ of characteristic $p$ and cardinality $q$. Let $\infty$ be a fixed place of $F$ and denote by $A\subset F$ the subring of functions, regular outside $\infty$. Finally let $\Gamma\subset \GL_2(F)$ be a congruence subgroup. In this situation, the work \cite{BGP} introduces a filtration on the space $S_{k,l}(\Gamma)$ of Drinfeld cusp forms of weight $k\ge2$, type $l\in\BZ$ and level $\Gamma$, that is preserved by the natural Hecke action. More precisely, we give such a filtration on the space of $\Gamma$-invariant harmonic cocycles $C_\har(V_{k,l})^\Gamma$ with values in a certain $\GL_2(F)$-representation $V_{k,l}$ that is finite-dimensional as an $F$-vector space, and we use the Hecke-equivariant residue isomorphism $S_{k,l}(\Gamma)\cong C_\har(V_{k,l})^\Gamma\otimes_F\BC_\infty$ of Teitelbaum, as in \cite[Theorem 1.7]{BGP} or \cite{Teitelbaum-Poisson}, to transfer the latter filtration to $S_{k,l}(\Gamma)$; here $\BC_\infty$ is the topological closure of an algebraic closure of the completion $F_\infty$ of $F$ at $\infty$. To be fully correct, one should develop the above setting adelically. But as later we shall only be concerned with $A=\BF_q[t]$, we refer to \cite{BGP} for details.

The filtration we constructed arouse from a filtration of the $\GL_2(F)$-representation $V_{k,l}$. The equi-characteristic representation theory of $\GL_2(F)$ is well-understood, e.g.~\cite[Section~10]{Bonnafe}. For each $k,l$, there is a unique Jordan-H\"older factor $L_{k,l}$ of $V_{k,l}$, its highest weight constituent, that is intrinsic to weight $k$, in the sense that all other subfactor are of the form $L_{k',l'}$ for some $k'<k$. So a natural question that arose in \cite{BGP} was to find a formula for $\dim_F  C_\har(L_{k,l})^\Gamma$. We note that for $\Gamma\subset\SL_2(A)$, the representations $V_{k,l}$ and $L_{k,l}$ do not depend on $l$, and so in that case, and in what follows, we omit it from the notation. 

If $\Gamma$ is $p'$-torsion free, i.e., if all elements of finite order have order a power of $p$, then from \cite{Teitelbaum-Poisson} one deduces $\dim_F  C_\har(L_{k})^\Gamma=c_{\Gamma}\dim_F L_{k}$ where $c_\Gamma\in\BN$ is a constant that depends only on $\Gamma$. Moreover $\dim_F L_{k}$ can be given explicitly in terms of the base $q$ digit expansion of $k$. However for large groups $\Gamma$, such as $\SL_2(A)$, in \cite[Subsection 4.2]{BGP} closed formulas were only found in the case $A=\BF_q[t]$ and for $q=2,3,5$. Numerical computations suggested an asymptotic formula
\begin{equation}\label{eq:asympt}
    \dim_F  C_\har(L_{k})^{\SL_2(A)}\approx\frac{\gcd(2,q^2-1)}{q^2-1} \dim_F L_{k}
\end{equation}
when $\dim_F L_{k}$ is `large', and for $q=2,3,5$ the formulas could be written as the floor function of a rational linear function in $\dim_F L_{k}$. Since this was not made explicit in \cite{BGP}, let us state this explicitly:
\begin{ThmInt}[{cf.~\autoref{p5}, \autoref{p3} and \autoref{p2}}]\label{prop:q=235}
    Let $A=\BF_q[t]$ and write $d_{k,q}:=\dim_F  C_\har(L_{k})^{\SL_2(A)}$. Then 
    \begin{enumerate}
        \item If $q=2$, then $\dim_F L_k\equiv \pm1\pmod{3}$ and
        $d_{k,q}= \big\lfloor \frac{\dim_F L_k+1}{3}\big\rfloor.$
        \item If $q=3$, then $d_{k,q}= \big\lfloor \frac{\dim_F L_k+1}{4}\big\rfloor$, if $k$ is even, and $d_{k,q}= 0$, if $k$ is odd, and for even $k$ one has $\dim_F L_k\pmod 4\in\{0,1,3\}$.
        \item If $q=5$, then $d_{k,q}= \big\lfloor \frac{\dim_F L_k+7}{12}\big\rfloor$, if $k$ is even, and $d_{k,q}= 0$, if $k$ is odd, and for even $k$ one has $\dim_F L_k\pmod {12}\in\{0,1,3,4,5,8,9\}$.
    \end{enumerate}
\end{ThmInt}
In this short note, we shall derive, for $A=\BF_q[t]$, a general formula for $d_{k,q}$ in terms of representation theoretic data from the tori in $\SL_2(A)$ and the base $q$-digits of $k$; see \autoref{Lk} and \autoref{Dasding} The resulting formulas are a lot more complicated than those above. We work out the cases $q=4,7,11$ in Subsection \ref{second_cases}. Building on \autoref{Dasding}, we also prove the asymptotic formula \eqref{eq:asympt} in \autoref{prop:asympt}. 

As suggested in \cite[Remark~4.13]{BGP}, here we use Brauer characters to replace a very special method built on \cite{Reduzzi}, and we completely rework and generalize the calculations of \cite[Section~4.2]{BGP}. Let us also note that we had some hopes that for arbitrary $q$ the dimension formulas we compute might take a shape as in the above theorem. However the formulas in Subsection \ref{second_cases} appear to be a lot more complicated, and we now doubt that there can be simple closed formulas for these dimensions; cf.~Remark~\ref{rem:no-simple-formula}. Finally, we want to mention that some of the work presented here is part of the bachelor thesis of the third named author.

{\bf Acknowledgements:} G.B. received funding by the Deutsche Forschungsgemeinschaft (DFG, German Research Foundation) TRR 326 \textit{Geometry and Arithmetic of Uniformized Structures}, project number 444845124. P.G. received funding by the Deutsche Forschungsgemeinschaft (DFG, German Research Foundation) -- project number 546550122.

\section{Recollections on splitting fields and Brauer characters}
Let $p$ be a prime number and let $G$ be a finite group.

\medskip

\begin{Def}

Let $A$ be a finite-dimensional $K$-algebra for a field $K$. 
Call an overfield $E\supset K$ a splitting field for $A$ if every simple $A_E:=A\otimes_KE$-module is absolutely simple.

Call a field $E$ with prime field $E_0$ a \emph{splitting field} for $G$ if $E$ is a splitting field for $E_0[G]$, i.e., if every simple $E$-representation of $G$ is absolutely irreducible.

\end{Def}
Let $m= \exp(G)$ be the exponent of $G$. For a field $E$, if $p=\Char E >0$, let $m'$ be the prime-to-$p$ part of $m$, and if $\Char E=0$, let $m'=m$.
\begin{Prop}[{\cite[Theorem 17.1]{CR81}}]
If the field $E$ contains a primitive $m'$-th root of unity, then $E$ is a splitting field for $G$.
\end{Prop}
In particular the minimal splitting fields of $G$ are finite field or number fields, that if needed can be chosen as subfields of a cyclotomic field.

One has the following elementary criterion for being a splitting field.
\begin{Prop}[{\cite[Proposition 7.15]{CR81}}]
Let $K$ be a field and let $A$ be a finite-dimensional $K$-algebra. 
Then for an extension field $E$ of $K$, the following statements are equivalent:
\begin{enumerate}
\item $E$ is a splitting field for $A$.
\item For every simple $A_E$-module M, the natural map $E\to \End_{A_E}(M)$ is an isomorphism 
\end{enumerate}
\end{Prop}

An element $g\in G$ is called $p$-regular if its order is prime to $p$. As $p$-regularity is preserved under conjugation, it is a property of conjugacy classes. 

\begin{Prop}[{\cite[Corollary 17.11]{CR81}}]
Let $E$ be a field of characteristic $p$ that is a splitting field for $G$. 
Then the number of isomorphism classes of absolutely simple $E[G]$-modules, is equal to the number of $p$-regular conjugacy classes of $G$.
\end{Prop}

Let now $\kappa\supset\BF_p$ be a finite field that is a splitting field for $G$. Let also $G^\preg\subset G$ be the union of the $p$-regular conjugacy classes of~$G$. 

To define Brauer characters of irreducible $\kappa[G]$-representation, let $\alpha\in\kappa^\times$ be a generator, and fix a chosen root of unity $\xi\in\BC^\times$ of exact order $\#\kappa^\times$. Then there is a unique homomorphism  $\kappa\to  \BC$ of multiplicative monoids that maps $\alpha$ to $\xi$, and hence $\alpha^m$ to $\xi^m$ and $0$ to $0$.

Let $V$ be an irreducible  $\kappa[G]$-representation. Then on $g\in G^\preg$ with eigenvalue sequence
$(\alpha^{m_i(g)})_{i=1,\ldots,\dim_\kappa V}$ for its action on $V$, with $m_i(g)\in\BZ$, define the Brauer character $\psi_V$ by 
\[\psi_V : g \mapsto \sum_{i=1,\ldots,\dim_\kappa V} \xi^{m_i(g)}.\] 
Setting for $\phi,\psi:G^\preg\to\BC$
\begin{equation}\label{eq-Hermitian}
    \langle\phi,\psi\rangle :=\langle\phi,\psi\rangle_{G,p} := \frac1{\#G}\sum_{g\in G^\preg}\phi(g)\overline\psi(g),
\end{equation}
defines a hermitian form on the space of maps $G^\preg\to\BC$.

The key result for us on Brauer characters is 
\begin{Thm}[{\cite[Proposition 10.2.1]{Webb}}]
\label{dimform}
Let $P$ and $U$ be irreducible $\kappa [G]$-modules with Brauer characters $\psi_P$ and $\psi_U$, respectively. 
If $P$ is projective, then one has
\[ \dim_\kappa \Hom_{\kappa[G]}(P,U) =\langle \psi_P,\psi_U\rangle.\]
\end{Thm}

\section{Conjugacy classes of the group \texorpdfstring{$\SL_2(\BF_q)$}{SLTwoFq}}
For the remainder of this article, let $G=\SL_2(\BF_q)$ where $q=p^e$ for some $e>0$.

Using the generalized Jordan canonical form for elements in $G$, and taking into account that in some cases $\GL_2(\BF_q)$ conjugacy are larger than $\SL_2(\BF_q)$-conjugacy classes, one obtains the following result. 
\begin{Prop}\label{G-conjclasses}
The conjugacy classes of $G$ have the following representatives:
\begin{enumerate}
\item The diagonal matrices $\begin{pmatrix} 1&0\\0&1\end{pmatrix}$ and $\begin{pmatrix} -1&0\\0&-1\end{pmatrix}$.    
\item The companion matrices $\begin{pmatrix} 0&-1\\1&b\end{pmatrix}$ of the polynomials $T^2-bT+1$ with $b\in\BF_q\setminus\{\pm2\}$. 
\item The companion matrices $\begin{pmatrix} 1&1\\0&1\end{pmatrix}$ and $\begin{pmatrix} -1&-1\\0&-1\end{pmatrix}$ of $(T-1)^2$ and of $(T+1)^2$.
\item The matrices from (c) multiplied on the right by $\begin{pmatrix} 1&z\\0&1\end{pmatrix}$ for a fixed $z\in\BF_q\setminus\BF_q^{\times 2}$.
\end{enumerate}
The $p$-regular classes are represented by the elements in (a) and (b). There are $q$ of them.
\end{Prop}
\begin{Rem} 
\begin{enumerate}
\item The matrices in \autoref{G-conjclasses} (a) and (b)  represent all semisimple elements of $G$.
\item The matrices in (b) can be decomposed into $2$ classes: (i) Companion matrices for reducible polynomials. (ii) Companion matrices for irreducible polynomials. 

The matrices belonging to (i) have eigenvalues $\{\zeta,\zeta^{-1}\}$ with $\zeta\in \BF_q^\times\setminus\{\pm1\}$. 
The eigenvalues of the matrices belonging to (ii) are pairs $\{\zeta,\zeta^q\}$ with $\zeta \in \BF_{q^2}\setminus\BF_q$ a $q+1$-th root of unity.
\end{enumerate}
\end{Rem}

\section{Equal characteristic representations of the group \texorpdfstring{$\SL_2(\BF_q)$}{SLTwoFq}}
Let $F\supset \BF_q$ be any, possibly infinite, overfield. 

\begin{Def}
Let $F[X,Y]$ be the space of polynomials in two variables $X$ and $Y$ over $F$. For $k\in \BN_0$ we let $\Delta_k:=F[X,Y]_k\subset F[X,Y]$ be the subspace of homogeneous polynomials of degree~$k$. The group $\GL_2(F)$ acts on $F[X,Y]$ via
\begin{equation*}
\begin{pmatrix}
a & b \\ c & d 
\end{pmatrix} \cdot X^iY^j = (aX+cY)^i(bX+dY)^j\quad\hbox{for \ }\begin{pmatrix}
a & b \\ c & d 
\end{pmatrix}\in\GL_2(F). 
\end{equation*}
\begin{enumerate}
\item The subspace $\Delta_k$ is a $\GL_2(F)$-subrepresentation. For $F=\BF_q$, we denote its restriction to $G$ by~$\overline{\Delta}_k$.
\item The $F$-span of $\{X^iY^{k-i}\,|\,\binom{k}{i}\not\equiv 0 \mod p\}$ is a $\GL_2(F)$-subrepresentation of $F[X,Y]_k$. We denote it by $L_k$. For $F=\BF_q$, its restriction to $G$ is denoted by $\overline{L}_k$.
\item We define the \textit{Steinberg module} $\overline{\text{st}}$ to be $\overline{L}_{q-1}$.
\end{enumerate}
\end{Def}

\begin{Def}
Consider the group homomorphism
\begin{equation*}
\tau : \operatorname{GL}_2(F)\to\operatorname{GL}_2(F) \,\,,\,\, \begin{pmatrix}
a & b \\ c&d
\end{pmatrix}\mapsto \begin{pmatrix}
a^p & b^p \\ c^p&d^p
\end{pmatrix}.
\end{equation*}
For an $F[\GL_2(F)]$-module $V$ and $s\ge0$, we define its $s$-fold Frobenius twist to be the $F[\operatorname{GL}_2(F)]$-module $V^{(s)}$ with underlying vector space $V$ and $\GL_2(F)$-action given by:
\begin{equation*}
\operatorname{GL}_2(F)\times V\to V\,\,,\,\, (g,v)\mapsto \tau^s(g)\cdot v.
\end{equation*}
\end{Def}
Observe for later that for Brauer characters one has
\begin{equation}\label{eq:CharOfFrob}
    \psi_{V^{(s)}}(g)=\psi_V(g^{p^s}).
\end{equation}
The next proposition sums up the most important properties of the modules defined above.
\begin{Prop} 
\label{simplemods}
Let $k$ be a non negative integer and $k = k_0 + k_1 p+ \dots + k_rp^r$ be a base $p$ expansion of $k$, meaning $0\leq k_i \leq p-1$ for all $i$ and $k_r\neq 0$. Then the following holds:
\begin{enumerate}
\item The modules $L_i$, $0\le i< \#F$, form a complete set of representatives of irreducible $F[\operatorname{SL}_2(F)]$-modules $V$ with $\dim_FV<\infty$. 
The $L_i$ are absolutely irreducible.
\item The modules $\overline{L}_0,\dots,\overline{L}_{q-1}$ form a complete set of representatives of irreducible $\BF_q[G]$-modules; they are absolutely irreducible and in particular, $\BF_q$ is a splitting field of $G$.
\item We have $\Delta_a=L_a$ precisely for $a=0,\ldots,p-1$, and 
\begin{equation*}
\bigotimes_{i=0}^r\Delta_{k_i}^{(i)}\to{\Delta}_{k}\,\,,\,\,f_0(X,Y)\otimes\dots\otimes f_r(X,Y)\mapsto \prod_{i=0}^r f_i(X^{p^i},Y^{p^i})
\end{equation*}
is a monomorphism of $\GL_2(F)$-representations which has image $L_k$.
\item The Steinberg module $\overline{\operatorname{st}}$ is a projective $\BF_q[G]$-module.
\end{enumerate}
\end{Prop}

\begin{proof}
See 
\cite[Lemma 8]{Pellarin-CharPReps} for (a), \cite[Proposition 10.2.4]{Bonnafe} for (c) and \cite[Proposition 10.1.8]{Bonnafe} for (d), and observe that (b) is a special case of (a), and that we may use $F=\BF_q$ to deduce the splitting field property.
\end{proof}
\begin{Rem}
    Regard $F^2$ as a space of column vectors equipped with a left action of $\GL_2(F)$ by left multiplication with matrices. Then $\Delta_k\cong \Sym^k F^{\oplus 2}$. The representation $V_k$ mentioned in the introduction is then given by its dual, i.e., $V_k=\Hom_F(\Delta_k,F)$.
\end{Rem}

\section{Brauer characters of certain \texorpdfstring{$\BF_q[\SL_2(\BF_q)]$}{SLTwoFq-Reps}-representations}

We want to apply \autoref{dimform} to obtain dimension formulas for $\Hom_{\BF_q[\SL_2(\BF_q)]}(\overline{\operatorname{st}}, \overline{L}_k)$. By $\mu_{q-1}$ and $\mu_{q+1}$ we denote the subgroups of $\BC^\times$ of elements of order dividing $q-1$ or $q+1$, respectively, or their reductions to $\BF_q^\times$ or $\BF_{q^2}^\times$.

In the following proposition, whenever $p=2$, expressions like $\pm 1$ are supposed to represent just the one element~$1$. Denote by $\Im(z)$ the imaginary part of $z\in\BC$. 

\begin{Prop}
\label{Deltak}
Let $\zeta \in \BF_{q^2}^\times$ be either a primitive $(q+1)$-th or $(q-1)$-th root of unity and let $\overline{\zeta}$ in $\BC^\times$ be a lift of the same order. Let $A_\zeta\in G$ denote the matrix 
\begin{equation*}
A_\zeta = 
\begin{cases}
 \pm I_2 &\text{ if }\zeta =\pm 1\\
 \begin{pmatrix}
 0 & -1 \\ 1& \zeta + \zeta^{-1}
\end{pmatrix}  & \text{ otherwise }
\end{cases}
\end{equation*}
The Brauer character of $\overline\Delta_k$ is given by 
\begin{equation*}
\psi_{\overline\Delta_k}(A_\zeta)=
\begin{cases}
(\pm1)^{k}(k+1) & \text{ if } \zeta = \pm 1 \\
\frac{\Im(\overline\zeta^{k+1})}{\Im(\overline\zeta)}
& \text{ otherwise.}
\end{cases}  
\end{equation*}
For $\overline{\st}=\otimes_{i=0,\ldots,e-1} \overline{\Delta}_{p-1}^{(i)}$ we have
\begin{equation*}
\psi_{\overline{\st}}(A_\zeta)=
\begin{cases}
q	 	&\text{ if } \zeta = \pm 1\\
1 		&\text{ if } \zeta \in \BF_q^\times\setminus \{\pm 1\} \\
-1 		&\text{ if } \zeta \in \mu_{q+1}\setminus \{\pm 1\}
\end{cases}  
\end{equation*}
\end{Prop}

\begin{proof}

Let $M$ denote the matrix by which $A_\zeta$ acts on $\Delta_k$ with respect to the basis $X^iY^{k-i}$, $i = 0,\dots,k$. If $\zeta = \pm 1$, then
$$A_\zeta\cdot X^iY^{k-i} = (\zeta X)^i(\zeta Y)^{k-i} = \zeta^k X^iY^{k-i}.$$
This means that the eigenvalues of $M$ are all $\zeta^{k}$. Their lifts sum up to $\overline{\zeta}^{k}(k+1)$.

If $\zeta \neq \pm 1$ is a $(q-1)$-th root of unity, $A_\zeta$ is conjugate to the diagonal matrix with entries $\zeta$ and $\zeta^{-1}$. The action of this element of $G$ is represented by $\operatorname{diag}(\zeta^{-k},\zeta^{2-k},\zeta^{4-k},\dots,\zeta^{2k-k})$, because
\begin{equation*}
\begin{pmatrix}
\zeta & 0 \\ 0 & \zeta^{-1}
\end{pmatrix}\cdot X^iY^{k-i} = (\zeta X)^i(\zeta^{-1} Y)^{k-i} = \zeta^{2i-k} X^iY^{k-i}
\end{equation*}
Thus the Brauer character is given by
$$\psi_{\overline{\Delta}_k}(A_\zeta)=\sum_{i=0}^k\overline{\zeta}^{-k}\overline{\zeta}^{2i} = \zeta^{-k}\frac{\overline{\zeta}^{2(k+1)}-1}{\overline{\zeta}^{2}-1}=\frac{\overline{\zeta}^{k+1}-\overline{\zeta}^{-k-1}}{\overline{\zeta}-\overline{\zeta}^{-1}}=\frac{\Im(\zeta^{k+1})}{\Im(\zeta)}.$$

The case $\zeta\in\mu_{q+1}\setminus\{\pm1\}$ is analogous. One can work with $F=\BF_{q^2}$ over which the action of $A_\zeta$ is diagonalizable, and the argument is the same as for $\zeta\in\mu_{q-1}\setminus\{\pm1\}$. 

To evaluate the Brauer character of $\overline{\st}$, observe that for $a=p-1$ we have $(-1)^a=1$ for any prime $p$. 
Now \autoref{simplemods}, equation~\eqref{eq:CharOfFrob} and the already proved parts imply that
\[  \psi_{\overline{L}_{q-1}}(A_\zeta)=\prod_{i=0}^{e-1} \psi_{\overline{\Delta}_{p-1}}(A_{\zeta^{p^i}})= 
\begin{cases}
\prod_{i=0}^{e-1} p	 	&\text{ if } \zeta = \pm 1,\\
\prod_{i=0}^{e-1} \frac{\Im(\overline\zeta^{p^{i+1}})}{\Im(\overline\zeta^{p^i})} 		&\text{ if } \zeta \neq \pm1.
\end{cases}  
\]
Now in the second case, one uses that $\zeta^{p^e}=\zeta$ for $\zeta\in\mu_{q-1}$ and $\zeta^{p^e}=\zeta^{-1}$ for $\zeta\in\mu_{q+1}$ and $\Im(\zeta^{-1})=-\Im(\zeta)$ to conclude. 
\end{proof}

From now on we will fix a non negative integer $k$.
Let $k = k_0+pk_1+\dots+p^rk_r$ be a base $p$ expansion of $k$ (with $k_r\neq 0$) and let $k = l_0+ql_1+\dots+q^tl_t$ be a $q$-expansion of $k$. For the sake of notation we put $k_i = 0$ for all $i > r$. In this situation we have $l_i = k_{ie}+pk_{ie+1}+\dots+p^{e-1}k_{ie+(e-1)}$. 
For $i\in \{0,\dots,p-1\}$ and $j\in\{0,\dots,e-1\}$ we define 
\begin{equation}\label{eq:mij-def}
m_{i,j} := \#\{m\in\{1,\dots,t\}|\,k_{j+me}=i\}.    
\end{equation}
This is the number of digits in the base $q$ expansion of $k$, which have an $i$ as the $j$-th digit in their base $p$ expansion. In the case where $p = q$ (where $j$ is always zero) $m_{i,0}$ is just the number of times $i$ occurs in the base $p$ expansion of $k$.

\begin{Prop}
\label{Lk}
With the notation from \autoref{Deltak} we get for the Brauer character $\psi_{\overline{L}_k}$ of $\overline{L}_k$:
$$\psi_{\overline{L}_k}(A_\zeta)= \prod_{i=0}^{p-1}\prod_{j=0}^{e-1}(\psi_{\overline{\Delta}_i}(A_{\zeta^{p^j}}))^{m_{i,j}}=
\begin{cases}
\prod_{i=0}^{p-1}\prod_{j=0}^{e-1}((\pm1)^{i}(i+1))^{m_{i,j}}
& \text{ if } \zeta = \pm 1 \\
\prod_{i=0}^{p-1}\prod_{j=0}^{e-1}
(\frac{\Im(\overline\zeta^{(i+1)p^j})}{\Im(\overline\zeta^{p^j})})^{m_{i,j}}
& \text{ otherwise.}
\end{cases}  
$$
\end{Prop}

\begin{proof}
By \autoref{simplemods} we have 
$$\overline{L}_k= \bigotimes_{i=0}^r\overline{\Delta}_{k_i}^{(i)}.$$ Since for any $\BF_q[G]$-module $V$, the $e$-fold Frobenius twist $V^{(e)}$ is equal to the module $V$ itself, the resulting tensor product can be simplified in terms of the $m_{i,j}$:
$$\overline{L}_k\cong \bigotimes_{i=0}^r{\overline{\Delta}}_{k_i}^{(i)}\cong \bigotimes_{i=0}^{p-1}\bigotimes_{j=0}^{e-1} (\overline{\Delta}_i^{(j)})^{\otimes m_{i,j}}$$
This translates to an equality of Brauer characters
$$\psi_{\overline{L}_k}(A_\zeta)= \prod_{i=0}^{p-1}\prod_{j=0}^{e-1}(\psi_{\overline{\Delta}_i^{(j)}}(A_{\zeta}))^{m_{i,j}},$$
and we conclude using formula~\eqref{eq:CharOfFrob} and \autoref{Deltak}.
\end{proof}

\begin{Thm}
\label{Dasding}
Let $\overline{V}$ be any finite dimensional $\BF_q[G]$-module with Brauer character $\psi_{\overline V}$. Then the multiplicity of $\overline{\st}$ in $\overline{V}$ is given by the following formula:
\begin{equation*}
\dim_{\BF_q}\Hom_{\BF_q[G]}(\overline{\st}, \overline{V})=\frac{1}{2}\left(\frac{1}{q-1}\sum_{\zeta \in \mu_{q-1}}\psi_{\overline{V}}(A_\zeta)-\frac{1}{q+1}\sum_{\zeta \in \mu_{q+1}}\psi_{\overline{V}}(A_\zeta)\right)
\end{equation*}
\end{Thm}

\begin{proof}
We write $G^\preg/\!\sim\,$ for the set of regular conjugacy classes ${}^Gg\subset G^\preg$. 
By \autoref{dimform} and \autoref{Deltak} regardless of the parity of $p$, we have:
\begin{eqnarray*}
\lefteqn{\dim_{\BF_q}\Hom_{\BF_q[G]}(\overline{\st}, \overline{V}) }\\
			&=&\frac{1}{\#G}\sum_{g\in G^\preg}\overline{\psi_{\overline{\st}}(g)}\psi_{\overline{V}}(g) \ = \ \frac{1}{\#G}\sum_{{}^Gg\in G^\preg/\!\sim}\#({}^Gg)\overline{\psi_{\overline{\st}}(g)}\psi_{\overline{V}}(g)\\
			&=&{\frac{1}{q(q^2-1)}}\left( q \sum_{\zeta\in \{\pm 1\} }\psi_{\overline{V}}(A_\zeta) \right. +{
            \frac{q(q+1)}{2}} \!\!\!\sum_{\zeta\in \mu_{q-1}\setminus\{\pm 1\}}\psi_{\overline{V}}(A_\zeta)-\left. \frac{q(q-1)}{2} \!\!\!\sum_{\zeta\in \mu_{q+1}\setminus\{\pm 1	\}}\psi_{\overline{V}}(A_\zeta)\right)\\
			&=&\frac{1}{2}\left(\frac{1}{q-1}\sum_{\zeta \in \mu_{q-1}}\psi_{\overline{V}}(A_\zeta)-\frac{1}{q+1}\sum_{\zeta \in \mu_{q+1}}\psi_{\overline{V}}(A_\zeta)\right)
\end{eqnarray*} 
\end{proof}

\begin{Rem}
With the notation from formula \eqref{eq-Hermitian}, the formula above can also be written as 
$$\dim_{\BF_q}\Hom_{\BF_q[G]}(\overline{\st}, \overline{V}) = \frac{1}{2}(\langle \psi_{\overline{\st}},\psi_{\overline{V}} \rangle_{\mu_{q-1},p} +\langle \psi_{\overline{\st}},\psi_{\overline{V}} \rangle_{\mu_{q+1},p})-(\psi_{\overline{V}}(I_2)+\psi_{\overline{V}}(-I_2))$$
if $p$ is odd and, if $p=2$, as
$$\dim_{\BF_q} \Hom_{\BF_q[G]}(\overline{\st}, \overline{V}) = \frac{1}{2}(\langle \psi_{\overline{\st}},\psi_{\overline{V}} \rangle_{\mu_{q-1},p} +\langle \psi_{\overline{\st}},\psi_{\overline{V}} \rangle_{\mu_{q+1},p})-\psi_{\overline{V}}(I_2).$$
\end{Rem}

Recall that as indicated in the introduction and following \cite{BGP} we were originally interested in the quantity
\[
d_{k,q}=\dim_F  C_\har(L_{k})^{\SL_2(A)}.
\]
By \cite[Proposition 4.1]{BGP}, we have
\[
d_{k,q}=\dim_F\Hom_{\BF_q[G]}(\overline{\st},L_k)=\dim_{\BF_q}\Hom_{\BF_q[G]}(\overline{\st},\overline{L}_k).
\]

Hence, we can use \autoref{Dasding} together with \autoref{Lk} to compute $d_{k,q}$. We will do precisely that for small values of $q$ in the next section.

\section{Computations for small \texorpdfstring{$q$}{q}}

\subsection{The cases \texorpdfstring{$q=2,3,5$}{q=2,3,5}}
\label{first_cases}
We first consider the cases $p=q=2,3,5$ that were studied in \cite{BGP} by different methods, see \cite[Proposition 4.9, Proposition 4.10 and Proposition 4.11]{BGP} for details. Only the calculation for $p = 5$ will be done in complete detail. The other results are obtained analogously. We observe that in the case where $p=q$ the second index of $m_{i,j}$ as in (\ref{eq:mij-def}) is always zero -- hence in this case we omit it from the notation and we write $m_{i,0}=m_i$. In the sequel, the expression $0^n$ for any natural number $n$ is equal to $1$ if $n=0$ and zero otherwise.

\begin{Lem}\label{lem:dimLk}
We have 
$$\dim_F L_k = \prod_{i=1}^{p-1}\prod_{j=0}^{e-1}(i+1)^{m_{i,j}}.$$
In particular for $q=p$ we obtain
$$\dim_F L_k = 2^{m_1}\cdot\dots\cdot p^{m_{p-1}}.$$
\end{Lem}

\begin{proof}
This is a direct consequence of \autoref{Lk}, since 
$\dim_F L_k =\dim_{\BF_q}\overline{L}_k=\psi_{\overline{L}_k}(I_2)$.
\end{proof}

\begin{Prop}[$q=5$]
\label{p5}
For $k$ even we have 
\[
d_{k,5}= \left\lfloor \frac{\dim_F L_k+7}{12}\right\rfloor
\]
and $\dim_F L_k\pmod {12}\in\{0,1,3,4,5,8,9\}$, and $d_{k,5}= 0$ if $k$ is odd.
\end{Prop}

\begin{proof}
Let $\zeta_{4}$ (respectively $\zeta_{6}$) be a primitive $4$th (respectively $6$th) root of unity in $\BF_{25}$. They lift to primitive roots of unity $\overline{\zeta}_4$ and $\overline{\zeta}_6$ in $\BC$. Without loss of generality we may choose these primitive roots of unity to be $\zeta_k = e^{\frac{2\pi i}{k}}$. 
\autoref{Deltak} allows us to calculate the character of $\Delta_k$ for $k= 1,\dots, 4$. The following table shows the results. Note that we only need to consider the elements $A_1$, $A_{-1}$, $A_{\zeta_4}$, $A_{\zeta_6}$, $A_{\zeta_6^2}$, because they represent the $p$-regular conjugacy classes.

{\renewcommand{\arraystretch}{1.2}
\begin{center}
\begin{tabular}{|c|c|c|c|c|}
\hline
$\zeta$ & $\psi_{\overline{\Delta}_1}(A_\zeta)$ & $\psi_{\overline{\Delta}_2}(A_\zeta)$ & $\psi_{\overline{\Delta}_3}(A_\zeta)$ & $\psi_{\overline{\Delta}_4}(A_\zeta)$\\
\hline
$1$&$2$&$3$&$4$&$5$\\
\hline
$-1$&$-2$&$3$&$-4$&$5$\\
\hline
$\zeta_4^{\pm 1}$&$0$&$-1$&$0$&$1$\\
\hline
$\zeta_6^{\pm 1}$&$1$&$0$&$-1$&$-1$\\
\hline
$(\zeta_6^2)^{\pm 1}$&$-1$&$0$&$1$&$-1$\\
\hline
\end{tabular}
\end{center}
}

For an arbitrary $k$ we may use \autoref{Lk} to compute the Brauer character of $L_k$. Applying \autoref{Dasding} yields
\begin{equation*}
\begin{split}
d_{k,5}=& 
													\frac{1}{8} (2^{m_1}3^{m_2}4^{m_3}5^{m_4}+(-2)^{m_1}3^{m_2}(-4)^{m_3}5^{m_4}+2(-1)^{m_2}0^{m_1+m_3})\\
										&-\frac{1}{12}(2^{m_1}3^{m_2}4^{m_3}5^{m_4}+(-2)^{m_1}3^{m_2}(-4)^{m_3}5^{m_4}\\
										&+2(-1)^{m_4}((-1)^{m_1}+(-1)^{m_3})).
\end{split}
\end{equation*}
Using \autoref{lem:dimLk} this simplifies to 
\begin{align}
\label{p5_form}
d_{k,5}= \frac{1}{24}((\dim_F L_k-4(-1)^{m_4+m_3}0^{m_2})(1+(-1)^{m_1+m_3})+6(-1)^{m_2}0^{m_1+m_3}).
\end{align}
Observe that since $5\equiv 1 \mod 2$, we have $k \equiv m_1+m_3\mod 2$ proving the result for odd $k$. Now let $k$ be even. Then so is $m_1+m_3$. $\dim_F L_k$ can never be congruent to $7$ or $11$ modulo $12$, because both seven and eleven are units modulo $12$, and $\dim_F L_k  = 2^{m_1}3^{m_2}4^{m_3}5^{m_4}$ can only be a unit if $m_1= m_2 = m_3 = 0$, in which case it is a power of $5$ and thus congruent to $1$ or $5$ modulo $12$. Considering $\dim_F L_k$ modulo $12$ gives:
$$\dim_F L_k \equiv 
\left\{\begin{array}{ll}
0 \mod 12 &\iff m_1+2m_3\geq 2 \text{ and } m_2 \geq 1\\
1 \mod 12 &\iff m_1 =m_2=m_3 = 0 \text{ and } m_4 \text{ is even}\\
2 \mod 12 &\iff m_1 = 1, m_2=m_3 = 0 \text{ and } m_4 \text{ is even}\\
3 \mod 12 &\iff m_1 = m_3 = 0 \text{ and } m_2 \text{ is odd} \\
4 \mod 12 &\iff m_2 = 0, m_1 +2m_3 \geq 2 \text{ and } m_1 + m_4 \text{ is even}  \\
5 \mod 12 &\iff m_1 =m_2=m_3 = 0 \text{ and }m_4 \text{ is odd}\\
6 \mod 12 &\iff m_1 = m_2 = 1 \text{ and } m_3 = 0\\
8 \mod 12 &\iff m_2 = 0, m_1 +2m_3 \geq 2 \text{ and } m_1 + m_4 \text{ is odd}  \\
9 \mod 12 &\iff m_1 = m_3 = 0 \text{ and } m_2 \text{ is even and }\neq 0\\
10 \mod 12 &\iff m_1 = 1, m_2=m_3 = 0 \text{ and } m_4 \text{ is odd}
\end{array}\right. $$
The cases $\dim_F L_k \equiv 2,6,10 \mod 12$ lead to contradictions, since $m_1+m_3$ must be even. Thus they do not occur. We obtain the assertion by using formula (\ref{p5_form}) in the remaining cases.
\end{proof}

\begin{Prop}[$q=3$]
\label{p3}
For $k$ even we have
\[
d_{k,3}= \left\lfloor \frac{\dim_F L_k+1}{4}\right\rfloor
\]
and $\dim_F L_k\pmod 4\in\{0,1,3\}$, and $d_{k,3}= 0$ if $k$ is odd.
\end{Prop}

\begin{proof}
An analogous computation as in the proof of \autoref{p5} yields
\begin{align}
\label{p3_form}
d_{k,3}=\frac{1}{8}((1+(-1)^{m_1})\dim_F L_k-2\cdot 0^{m_1}(-1)^{m_2})
\end{align}
Since $3\equiv 1 \mod 2$, we have $k \equiv m_1\mod 2$ and we deduce that $d_{k,3}=0$ if $k$ is odd. Let $k$ be even. Then so is $m_1$. From \autoref{lem:dimLk} we obtain the following information on $m_1$ and $m_2$:
$$\dim_F L_k \equiv 
\left\{\begin{array}{ll}
0 \mod 4 &\iff m_1 \geq 2\\
1 \mod 4 &\iff m_2 \text{ is even and } m_1 = 0\\
2 \mod 4 &\iff m_1 = 1\\
3 \mod 4 &\iff m_2 \text{ is odd and } m_1 = 0\\
\end{array}\right. $$
The case $\dim L_k \equiv 2 \mod 12$ leads to a contradiction, since $m_1$ must be even. Applying formula (\ref{p3_form}) with this information yields the desired result.
\end{proof}

\begin{Prop}[$q=2$]
\label{p2}
We have $\dim_F L_k\equiv \pm1\pmod{3}$ and
\[
d_{k,2}= \left\lfloor \frac{\dim_F L_k+1}{3}\right\rfloor.
\]
\end{Prop}

\begin{proof}
An analogous computation as in the proof of \autoref{p5} yields
\begin{equation*}
d_{k,2}=\frac{1}{3}(\dim_F L_k-(-1)^{m_1})
\end{equation*}
and the result follows upon observing that $\dim_F L_k = 2^{m_1} \equiv (-1)^{m_1} \mod 3$.
\end{proof}

\subsection{The cases \texorpdfstring{$q=4,7,11$}{q=4,7,11}}
\label{second_cases}
The following formula is obtained in complete analogy with the results of the previous subsection. 
\begin{Prop}[$q=7$]
\label{p7}
We have
\begin{align*}
d_{k,7}= &\frac{1}{48}((1+(-1)^{m_1+m_3+m_5})\dim_F L_k+8\cdot 0^{m_2+m_5}(-1)^{m_1+m_4}(1+(-1)^{m_1+m_3})\\
	&-6\cdot 0^{m_3}((-1)^{m_2+m_6}0^{m_1+m_5}+2^\frac{m_1+m_5}{2}(1+(-1)^{m_1+m_5})(-1)^{m_1+m_4+m_6}).
\end{align*}
\end{Prop}
In order to state the results for $q=11$ and $q=4$ recall that the sequence of \emph{Lucas numbers} $(F_n)_{n\in\mathbb{N}}$ is defined as
\[
F_n:=\left\{
\begin{aligned}
&2, & n=0,\\
&1, & n=1,\\
&F_{n-1}+F_{n-2}, &n\geq 2.
\end{aligned}
\right.
\]

\begin{Prop}[$q=11$]
\label{p11}
We have
\begin{align*}
d_{k,11}=&\frac{1}{120}(\dim_F L_k(1+(-1)^{m_1+m_3+m_5+m_7+m_9})\\
	&-10(((-1)^{m_1+m_4+m_7+m_{10}}+(-1)^{m_3+m_4+m_9+m_{10}})\cdot 0^{m_2+m_5+m_8})\\
	&+0^{m_1+m_3+m_5+m_7+m_9}(-1)^{m_2+m_6+m_{10}}\\
	&+0^{m_5}(-1)^{m_1+m_3+m_6+m_8+m_{10}}2^{m_2+m_8}\\
	&\cdot 3^{\frac{m_1+m_3+m_7+m_9}{2}}(1+(-1)^{m_1+m_3+m_7+m_9}))\\
	&+12\cdot 0^{m_4+m_9}(-1)^{m_8+m_6}\\
	&\cdot ((-1)^{m_3+m_1}+(-1)^{m_5+m_7})F_{m_1+m_2+m_6+m_7}).\\
\end{align*}
\end{Prop}

\begin{proof}
Performing the same calculations as in \autoref{p5} leads to the following values for the Brauer characters of $\Delta_i$:

\renewcommand{\arraystretch}{1.2}
\begin{center}
\begin{tabular}{|c|c|c|c|c|c|c|c|c|c|c|}
\hline
\diagbox%
[width=\dimexpr \textwidth/16+\tabcolsep\relax, height=1cm]{$\zeta$}{$i$} &$1$&$2$&$3$&$4$&$5$&$6$&$7$&$8$&$9$&$10$\\
\hline
$1$&$2$&$3$&$4$&$5$&$6$&$7$&$8$&$9$&$10$&$11$\\
\hline
$-1$&$-2$&$3$&$-4$&$5$&$-6$&$7$&$-8$&$9$&$-10$&$11$\\
\hline
$\zeta_{12}^{\pm 1}$&$\sqrt{3}$&$2$&$\sqrt{3}$&$1$&$0$&$-1$&$-\sqrt{3}$&$-2$&$-\sqrt{3}$&$-1$\\
\hline
$\zeta_{12}^{\pm 2}$&$1$&$0$&$-1$&$-1$&$0$&$1$&$1$&$0$&$-1$&$-1$\\
\hline
$\zeta_{12}^{\pm 3}$&$0$&$-1$&$0$&$1$&$0$&$-1$&$0$&$1$&$0$&$-1$\\
\hline
$\zeta_{12}^{\pm 4}$&$-1$&$0$&$1$&$-1$&$0$&$1$&$-1$&$0$&$1$&$-1$\\
\hline
$\zeta_{12}^{\pm 5}$&$-\sqrt{3}$&$2$&$-\sqrt{3}$&$1$&$0$&$-1$&$\sqrt{3}$&$-2$&$\sqrt{3}$&$-1$\\
\hline
$\zeta_{10}^{\pm 1}$&$\ro$&$\ro$&$1$&$0$&$-1$&$-\ro$&$-\ro$&$-1$&$0$&$1$\\
\hline
$\zeta_{10}^{\pm 2}$&$-\ri$&$\ri$&$-1$&$0$&$1$&$-\ri$&$\ri$&$-1$&$0$&$1$\\
\hline
$\zeta_{10}^{\pm 3}$&$\ri$&$\ri$&$1$&$0$&$-1$&$-\ri$&$-\ri$&$-1$&$0$&$1$\\
\hline
$\zeta_{10}^{\pm 4}$&$-\ro$&$\ro$&$-1$&$0$&$1$&$-\ro$&$\ro$&$-1$&$0$&$1$\\
\hline
\end{tabular}
\end{center}
Here $\varphi = \frac{1+\sqrt{5}}{2}$ is the golden ratio and $\widehat{\varphi} = -\varphi^{-1}$. Thus the formula for $d_{k,11}$ contains the term 
\begin{align*}
\frac{1}{10} \,0^{m_4+m_9}	&((-1)^{m_3+m_8}(\ri^{m_2+m_6}(-\ri)^{m_1+m_6}+(-\ro)^{m_1+m_6}\ro^{m_2+m_7})\\	
						&+(-1)^{m_5+m_8}(\ri^{m_1+m_2}(-\ri)^{m_6+m_7}+(-\ro)^{m_6+m_7}\ro^{m_1+m_2}).
\end{align*}
This simplifies to 
$$\frac{1}{10} \,0^{m_4+m_9}(-1)^{m_8+m_6} ((-1)^{m_3+m_1}+(-1)^{m_5+m_7})(\ri^{m_1+m_2+m_6+m_7}+\ro^{m_1+m_2+m_6+m_7}).$$

It is easily checked that $\ri^n+\ro^n$ satisfies the recursive condition of the Lucas numbers and that $\ri^0+\ro^0=2$ and $\ri+\ro=1$. Thus $\ri^n+\ro^n= F_n$ and the assertion follows.
\end{proof}
We end this section by stating the result for $q=4$ without proof.
\begin{Prop}[$q=4$]
\label{q4}
We have 
\begin{equation*}
    d_{k,4}= \frac{1}{15}(\dim_F L_k+5(-1)^{m_{1,0}+m_{1,1}}-3(-1)^{m_{1,0}}F_{m_{1,1}-m_{1,0}}).
\end{equation*}
\end{Prop}

\begin{Rem}\label{rem:no-simple-formula}
Contrary to the cases studied in the previous subsection, in general there seems to be no formula for $d_{k,q}$ involving only the quantity $\dim_F L_k$.
\end{Rem}

\section{Asymptotic behavior}
We begin with some preparations.
\begin{Lem}\label{lem:FirstBound}
    Let $l\in\{1,\ldots,p\}$ and let $\lambda_0\in(0,\frac{\pi}{2l}]$. 
    Then the following hold:
    \begin{enumerate}
        \item For $l$ fixed, one has 
        \[ \frac{\sin(l \lambda_0)}{l\sin(\lambda_0)}>\left| \frac{\sin(l \lambda)}{l\sin(\lambda)}\right|\hbox{ \ for all }\lambda\in
        (\lambda_0,{\textstyle\frac{\pi}2}].
        \]
        \item The function $\{1,\ldots,p\}\to\BR_{\ge0},l\mapsto \big|\frac{\sin(l\lambda_0)}{l\sin(\lambda_0)}\big|$ is strictly decreasing in $l$.
    \end{enumerate}
\end{Lem}
\begin{proof}
    To prove (a), consider the function $f:(\lambda_0,\frac{\pi}2]\to\BR,\lambda\mapsto \frac{\sin(l \lambda)}{l\sin(\lambda)}$. Suppose first that $\lambda_0<\lambda<\frac\pi{l}$. On this domain, we have 
    \[f'(\lambda)\cdot\sin^2(\lambda)=l\cos(l\lambda)\sin(\lambda)-\sin(l\lambda)\cos(\lambda)<l\cos(\lambda)\sin(\lambda)-l\sin(\lambda)\cos(\lambda)=0,\]
    because $\cos(l\lambda)<\cos(\lambda)$ since $\cos$ is strictly decreasing on $(0,\pi)$, and $\frac1l\sin(l\lambda)<\sin(\lambda)$ because $\sin$ is concave on $(0,\pi)$. This shows that (a) holds for $\lambda\in(\lambda_0,\frac\pi{l}]$. For $\frac\pi{l}\le\lambda\le\frac{\pi}2$, we choose $b\in\BN$ such that $0\le\lambda':=\lambda-b\frac\pi l\le\frac\pi l$, and let $\lambda''=\max \{ \lambda',\frac\pi l-\lambda'\}$, so that $|\sin(l\lambda)|=\sin(l\lambda'')$. Then 
    \[ \left| \frac{\sin(l \lambda)}{l\sin(\lambda)}\right|=\frac{\sin(l \lambda'')}{l\sin(\lambda'')}\frac{\sin(\lambda'')}{\sin(\lambda)}< \frac{\sin(l \lambda')}{l\sin(\lambda')}< \frac{\sin(l \lambda_0)}{l\sin(\lambda_0)},\]
    where the last inequality follows from the case $\lambda_0\le\lambda\le\frac\pi l$ already treated, because $\lambda_0\le \frac\pi{2l}$.

    Regarding (b) note that $l\lambda_0\in(0,\frac\pi2]$. Now $\sin$ is concave for on $[0,\pi]$. It follows that $l\mapsto \big|\frac{\sin(l\lambda_0)}{l\sin(\lambda_0)}\big|$ is strictly decreasing in $l$.
    \end{proof}

    The roots of unity $\ol\zeta$ that occur in the formula in \autoref{Dasding} are the powers of $\zeta_{q\pm1}$. We shall need certain bounds on $|\Im(\ol\zeta^{(l+1)p^j})|$ for $l=0,\ldots,p-1$ and $j=0,\ldots,e-1$ for such $\ol\zeta$. Observe first that when $\ol\zeta$ traverses all powers of $\zeta_{q\pm1}$, then so does $\ol\zeta^{p^j}$ and vice versa. Moreover $|\Im(\ol\zeta)|=|\Im(-\ol\zeta)|=|\Im(\ol\zeta^{-1})|=|\Im(-\ol\zeta^{-1})|$. One concludes that it suffices to consider the powers $\ol\zeta$ of $\zeta'_{q\pm1}=\zeta_{\lcm(2,q\pm1)}$ with $\arg(\ol\zeta)\in[0,\frac\pi2]$; here, for $z\in\BC^\times$ we denote by $\arg(z)\in [0,2\pi)$ the unique angle such that $z=|z|e^{i\arg(z)}$.
    \begin{Lem}\label{lem:SecondBound}
        Let $l\in\{2,\ldots,p\}$, suppose $q\pm1\ge4$, and let $\ol\zeta$ be a power of $\zeta_{q\pm1}$ different from $\pm1$. Define 
        $m=\max\{4p,\lcm(2,q\pm1)\}$. Then for any $i\in\BN_0$ one has
        \[ 
        \left|\frac{\Im(\ol\zeta^{lp^i})}{l\Im(\ol\zeta^{p^i})} \right|\le\cos(2\pi/m)
        =:C_{q\pm1}<1.
        \]
    \end{Lem}
    \begin{proof}
        As observed before the lemma, it suffices to bound $\left|\frac{\Im(\ol\zeta^l)}{l\Im(\ol\zeta)} \right|$ for $\ol\zeta$ a power of $\zeta'_{q\pm1}$ with $0<\lambda:=\arg(\ol\zeta)\le \frac\pi2$, and let $\lambda_0:=2\pi/m$. Then $\lambda_0\in(0,\frac\pi{2l}]$, $\lambda\in[\lambda_0,\frac\pi2]$ and 
        $$\left|\frac{\Im(\ol\zeta^l)}{l\Im(\ol\zeta)} \right|=\frac{\sin(l\lambda)}{l\sin(\lambda)}.$$
        It follows from \autoref{lem:FirstBound} that
        \[\frac{\sin(l\lambda)}{l\sin(\lambda)}\le \frac{\sin(l\lambda_0)}{l\sin(\lambda_0)}\le \frac{\sin(2\lambda_0)}{2\sin(\lambda_0)}=\cos(\lambda_0)=\cos(2\pi/m)<1,\]
        and this completes the proof.
     \end{proof}

The following result is an immediate consequence of \autoref{Lk}, \autoref{lem:dimLk} and \autoref{lem:SecondBound}.
    \begin{Lem}\label{lem:Asympt1}
    With the notation from the \autoref{Lk} we get for the Brauer character $\psi_{\overline{L}_k}$ of $\overline{L}_k$ and for $\zeta\neq\pm1$ a root of unity of order dividing $q\pm1$ the bound:
    
$$\left|\frac{\psi_{\overline{L}_k}(A_\zeta)}{\dim_F L_k}\right|= 
\prod_{i=0}^{p-1}\prod_{j=0}^{e-1}
\left|\frac{\Im(\overline\zeta^{(i+1)p^j})}{(i+1)\Im(\overline\zeta^{p^j})}\right|^{m_{i,j}}<C_{q\pm1}^{\sum_{i=1}^{p-1}\sum_{j=0}^{e-1}m_{i,j}}.
$$
\end{Lem}

    \begin{Lem}
        Suppose $p$ is an odd prime. Let $k\ge0$ and define the tuple $(m_{i,j})_{i=0,\ldots,p-1,j=0,\ldots,e-1}$ according to \eqref{eq:mij-def} Then one has 
        \[ k\equiv\sum_{i=0}^{p-1}i\sum_{j=0}^{e-1}m_{i,j} \equiv\sum_{i=0,i\,\textrm{odd}}^{p-1}\,\,\sum_{j=0}^{e-1}m_{i,j}  \pmod{2}.\]
    \end{Lem}
    \begin{proof}
        Recall that $m_{i,j}$ counts the number of occurrences of the digit $i$ in the base $p$ expansion of $k$ as the coefficient of a $p$-power of the form $p^{j+mp^e}$ for some $m\in\BN_0$. Let $K_{i,j}$ be the set of such $p^{j+mp^e}$. Using that $p^m\equiv1\pmod 2$ for any $m\in\BN_0$, it follows that modulo $2$ we have
        \[k=\sum_{i=0}^{p-1}\sum_{j=0}^{e-1} \sum_{b\in K_{i,j}} i p^b\equiv \sum_{i=0}^{p-1}i\sum_{j=0}^{e-1} \#K_{i,j}\equiv \sum_{i=0,i\,\text{odd}}^{p-1}\,\,\sum_{j=0}^{e-1} m_{i,j}.\]
    \end{proof}
\begin{Prop}
\label{prop:asympt}
    If $p$ and $k$ are odd, then $d_{k,q}=0$. If $p$ or $k$ are even, then for $\dim_F L_k\to\infty$, one has 
    \begin{align*}
    \frac{d_{k,q}}{\dim_F L_k} \to 
    \begin{cases}\frac1{q^2-1},&\textrm{if }p=2,\\
    \frac2{q^2-1},&\textrm{if $p$ is odd.}\end{cases}
    \end{align*}
\end{Prop}
    \begin{proof}
        Note first that from \autoref{lem:dimLk} it follows that $M_k:=\sum_{i=1}^{p-1}\sum_{j=0}^{e-1}m_{i,j}\to\infty $, when $\dim_F L_k\to\infty$, and hence in this case one has
        \[C_{q\pm1}^{M_l}\to0.\]
        It follows from \autoref{Dasding}, \autoref{Lk} and \autoref{lem:Asympt1} that
        \[ \frac{d_{k,q}}{\dim_F L_k}=\frac{1}{2}\left(\frac{1}{q-1}\sum_{\zeta \in \mu_{q-1}\cap\{\pm1\}}\psi_{\ol L_k}(A_\zeta)-\frac{1}{q+1}-\sum_{\zeta \in \mu_{q+1}\cap\{\pm1\}}\psi_{\ol L_k}(A_\zeta) \right)\]
        converges to 
        \[\frac{1}{2}\left(\frac{1}{q-1}\sum_{\zeta \in \mu_{q-1}\cap\{\pm1\}}1-\frac{1}{q+1}\sum_{\zeta \in \mu_{q+1}\cap\{\pm1\}}1 \right)\]
        for $\dim_F L_k\to\infty$. The inner sum is $\frac2{q-1}-\frac2{q+1}$ if $p$ is odd, and half of that if $p$ is even. The assertion of the proposition for $p$ even or for $k$ even is now immediate. 

        If on the other hand, $p$ and $k$ are odd, then one verifies that, one the one hand, the expression on the left is invariant under replacing $\zeta$ by $-\zeta$, but that on the other and, this change also produces a sign change. It follows that then the expression on the left is zero, as had to be shown. 
    \end{proof}

\printbibliography
\end{document}